23-Sep-15

# Jacobi's epsilon and zeta function for moduli outside the interval [0, 1]

Milan Batista

milan.batista@fpp.uni-lj.si

(Jun-Aug 2015)

**Abstract**

The formulas that relate Jacobi's Epsilon and Zeta function with real moduli in the interval $(1,\infty)$ or with pure imaginary moduli to elliptic functions with moduli in the interval $[0,1]$ are derived.

**1 Introduction**

The NIST handbook [1; 2] provide formulas that relate elliptic integrals and Jacobian elliptical functions *sn*, *cn*, *dn* with real moduli in the interval $(1,\infty)$ or with pure imaginary moduli to functions with moduli in the interval $[0,1]$. The aim of this paper is derivation of such relations also for Jacobi's Epsilon function and Jacobi's Zeta function.

The Jacobi's epsilon function is defined by [2; 3].

$$\varepsilon(x,k) \equiv E(\text{am}(x,k),k) \equiv \int_0^x \text{dn}^2(t,k)dt \tag{1}$$

where *x* is real variable, $k \in [0,1]$ is modulus, *E* is incomplete elliptic integral of second kind, *am* is Jacobi's amplitude function. Zeta function *Z* is periodic part of $\varepsilon$ and is defined by [2; 3].

$$Z(x,k) \equiv \varepsilon(x,k) - \frac{E(k)}{K(k)}x \tag{2}$$

Both above functions are odd with respect to *x*

$$\varepsilon(-x,k) = -\varepsilon(x,k), \quad Z(-x,k) = -Z(x,k) \tag{3}$$

Also, since

$$K(-k) = K(k), \quad E(-k) = E(k), \quad \text{dn}(x,-k) = \text{dn}(x,k) \tag{4}$$

we have

$$\varepsilon(x,-k) = \varepsilon(x,k), \quad Z(x,-k) = Z(x,k) \tag{5}$$

We now stipulate that (1) and (2) are valid for any real modulus $|k|>1$ and for any pure imaginary modulus, i.e., for modulus of the form i$k$, where $i^2 = -1$ and $k$ is real. Our task is to find appropriate formulas that will give these definitions a sense.

## 2 The case when |k| > 1

We begin with relation Jacobi reciprocal relation ( [2],Eq 22.17.14)

$$\mathrm{dn}(x,k) = \mathrm{cn}(kx,k^{-1})$$

and identity (see for instance [4], Eq 2.1.11)

$$\mathrm{dn}^2(x,k^{-1}) - k^{-2}\mathrm{cn}^2(x,k^{-1}) = 1 - k^{-2}$$

With these a straightforward computation yields

$$\int_0^x \mathrm{dn}^2(t,k)\,dt = \int_0^x \mathrm{cn}^2(kt,k^{-1})\,dt$$
$$= \frac{1}{k}\int_0^{kx} \mathrm{cn}^2(\tau,k^{-1})\,d\tau$$
$$= k\int_0^{kx} \mathrm{dn}^2(\tau,k^{-1})\,d\tau + (1-k^2)x$$

From this, on using (1) and (2), we get

$$\varepsilon(x,k) = k\varepsilon(kx,k^{-1}) + (1-k^2)x \quad (k>1) \tag{6}$$

or, using (2),

$$\varepsilon(x,k) = \left[k^2 \frac{E(k^{-1})}{K(k^{-1})} + 1 - k^2\right]x + kZ(kx,k^{-1}) \quad (k>1) \tag{7}$$

The last two formulas are required formulas by which we can compute $\varepsilon(x,k)$ when $k>1$.

To obtain te formula for Z when $k>1$ we start with (2). From the formulas 19.7.3 in [1] we, by replacing $k$ with $k^{-1}$, obtain when $k>1$

$$K(k) = k^{-1}\left[K(k^{-1}) \pm iK(k'^{-1})\right]$$

$$E(k) = k\left[E(k^{-1}) - k'^{-2}K(k^{-1})\right] \pm ik\left[E(k'^{-1}) - k^{-2}K(k'^{-1})\right]$$

where

$$k' \equiv \frac{k}{\sqrt{k^2-1}}$$

Hence

$$\frac{E(k)}{K(k)} = 1 + \frac{k^2\left[K(k^{-1})(E(k^{-1})-K(k^{-1}))-E(k'^{-1})K(k'^{-1})\right]}{K^2(k^{-1})+K^2(k'^{-1})}$$
$$\pm i\frac{k^2\left[K(k'^{-1})(E(k^{-1})-K(k^{-1}))+E(k'^{-1})K(k^{-1})\right]}{K^2(k^{-1})+K^2(k'^{-1})} \quad (k>1)$$

Combining this with (2) and (6) yields required formula

$$Z(x,k) = kZ(kx,k^{-1}) + \frac{k^2K^2(k'^{-1})}{K^2(k^{-1})+K^2(k'^{-1})}\left[\frac{E(k^{-1})}{K(k^{-1})} + \frac{E(k'^{-1})}{K(k'^{-1})} - 1\right]x$$
$$\pm i\frac{k^2K(k^{-1})K(k'^{-1})}{K^2(k^{-1})+K^2(k'^{-1})}\left[\frac{E(k^{-1})}{K(k^{-1})} + \frac{E(k'^{-1})}{K(k'^{-1})} - 1\right]x \quad (8)$$

Evidently Z is complex for $k>1$.

On replacing $k$ with $-k$ in (6) and (8) we obtain, using (4) and (5)For, we by using (6) and (8) we obtain

$$\varepsilon(x,-k) = \varepsilon(x,k), \quad Z(x,-k) = Z(x,k) \quad (k>1) \tag{9}$$

In this way we found that formulas (6), (7) and (8) are valid for any real $|k|>1$.

### 3 The case when *k* is pure imaginary number

We start with the following relation ([2], Eq 22.17.8)

$$\text{dn}(x,ik) = \frac{1}{\text{dn}(x/k_1',k_1)}$$

where

$$k_1 = \frac{k}{\sqrt{1+k^2}}, \quad k_1' = \frac{1}{\sqrt{1+k^2}} \tag{10}$$

Substituting this into (1) we calculate, using (2),

$$\varepsilon(x,ik) = \int_0^x \frac{1}{\text{dn}^2(t/k_1',k_1)} dt$$
$$= k_1' \int_0^{x/k_1'} \frac{d\tau}{\text{dn}^2(\tau,k_1)}$$
$$= \frac{k_1'E(\text{am}(x/k_1',k_1),k_1)}{1-k_1^2} - \frac{k_1'k_1^2}{1-k_1^2}\frac{\text{sn}(x/k_1',k_1)\text{cn}(x/k_1',k_1)}{\text{dn}(x/k_1',k_1)}$$
$$= \frac{E(k_1)}{k_1'^2 Z(k_1)}x + \frac{Z(x/k_1',k_1)}{k_1'} - \frac{k_1^2}{k_1'}\frac{\text{sn}(x/k_1',k_1)\text{cn}(x/k_1',k_1)}{\text{dn}(x/k_1',k_1)}$$

Using identity (see for instance [2])

$$Z(x+K,k) = Z(x,k) - k^2 \frac{\text{sn}(x,k)\text{cn}(x,k)}{\text{dn}(x,k)}$$

we obtain required formula

$$\varepsilon(x,ik) = \frac{E(k_1)}{k_1'^2 Z(k_1)} x + \frac{Z(x/k_1' + K(k_1), k_1)}{k_1'} \tag{11}$$

To extend $Z$ function for pure imaginary modulus we start with relations 19.7.2 in [1]

$$K(ik) = k_1' K(k_1), \quad E(ik) = \frac{E(k_1)}{k_1'}$$

where $k_1$ and $k_1'$ are given by (10). We have

$$\frac{E(ik)}{K(ik)} = \frac{E(k_1)}{k_1'^2 K(k_1)}$$

Combining this with (2) and (11) gives

$$Z(x,ik) = \frac{Z(x/k_1' + K(k_1), k_1)}{k_1'} \tag{12}$$

Note that (11) and (12) are valid for any $k$ since both $k_1$ and $k_1'$ are less than one.

### 4 Numerical comparisons

The present formulas were implemented into computer program where for calculation of elliptic integrals and $Z$ function we use Carlson algorithms [5]. In Tables 1,2,3,4 we give comparison of calculation of $\varepsilon(x,k)$ and $Z(x,k)$ by present method, Maple, MuPAD and numerical integration with Newton-Cotes method (Quanc8). We note that MuPAD does not use $k$ but parameter $m = k^2$. As can be seen from the tables results agrees perfectly.

**Table 1.** Comparison of calculation of $\varepsilon(x,k)$

| x | k | Present | Maple/MuPAD* | Quanc8 |
|---|---|---|---|---|
|  | 0.5 | 0.490203 | 0.490203 | 0.490203 |
| 0.5 | 1 | 0.462117 | 0.462117 | 0.462117 |
|  | 2 | 0.367975 | 0.367975 | 0.367975 |

**Table 2.** Comparison of calculation of $\varepsilon(x,ik)$

| x | k | Present | Maple/MuPAD | Quanc8 |
|---|---|---|---|---|
|  | 0.5 | 0.510020 | 0.510020 | 0.510020 |
| 0.5 | 1 | 0.541445 | 0.541445 | 0.541445 |
|  | 2 | 0.689051 | 0.689051 | 0.689051 |

**Table 3.** Comparison of calculation of $Z(x,k)$

| x   | k   | Present              | Maple/MuPAD*         |
|-----|-----|----------------------|----------------------|
|     | 0.5 | 0.054948             | 0.054948             |
| 0.5 | 1   | 0.462117             | 0.462117             |
|     | 2   | 0.663361 - 0.419309i | 0.663361 - 0.419309i |

**Table 4.** Comparison of calculation of $Z(x,ik)$

| x   | k   | Present   | Maple/MuPAD* |
|-----|-----|-----------|--------------|
|     | 0.5 | -0.050738 | -0.050738    |
| 0.5 | 1   | -0.187029 | -0.187029    |
|     | 2   | -0.616203 | -0.616203    |

## 5 Application

As simple application of present formulas consider Euler's flexural elastica which has parametric form given by Love [6] (Ch 19, Eq 12)

$$x(u) = \frac{1}{\omega}\left\{-u + 2\left[\varepsilon(u+K,k) - \varepsilon(K)\right]\right\}, \quad y(u) = -\frac{2k}{\omega}\mathrm{cn}(u+K,k) \quad (k<1) \tag{13}$$

where $u = \omega s$ and $\omega$ is a parameter. To obtain in-flexural elastica Love integrate modified differential equations of elastica. However we can obtain in-flexural elastic from (13) by using (6). We have

$$x = \frac{1}{\omega k}\left\{(1-2k^2)ku + 2k^2\left[\varepsilon\left(ku+kK(k^{-1}),k^{-1}\right) - \varepsilon\left(kK(k^{-1}),k^{-1}\right)\right]\right\}, \quad y = -\frac{2k}{\omega}\mathrm{dn}\left(ku+kK(k^{-1}),k^{-1}\right)$$

To obtain symmetric curve we transform $u + K(k^{-1}) \to u$ and $x + \frac{2k^2}{\omega k}\varepsilon\left(kK(k^{-1}),k^{-1}\right) \to x$. In this way we obtain Love formulas (Ch 19, Eq 16)

$$x(u) = \frac{1}{\omega k}\left\{(1-2k^2)ku + 2k^2\varepsilon(ku,k^{-1})\right\}, \quad y(u) = -\frac{2k}{\omega}\mathrm{dn}(ku,k^{-1}) \quad (k>1) \tag{14}$$

23-Sep-15